\documentclass[preprint, 12pt]{elsarticle}

\journal{Journal of Combinatorial Theory, Series A}

\usepackage{isabelle,isabellesym}
\isabellestyle{tt}

\usepackage{etoolbox}
\newtoggle{wip}
\togglefalse{wip}

\iftoggle{wip}{
\usepackage{todonotes}

\usepackage{xcolor}

	\usepackage[mark,dirty={Modified!}]{gitinfo2}

}
{
	\usepackage[disable]{todonotes}
}

\RequirePackage[l2tabu, orthodox]{nag}
\usepackage{microtype}

\usepackage[utf8]{inputenc}
\usepackage[english]{babel}
\usepackage{graphicx}

\usepackage{hyperref}
\usepackage{amsmath,amsthm,amssymb}
\usepackage{enumerate}
\usepackage{color}
\usepackage[nameinlink]{cleveref} 
\usepackage{tikz}

\usepackage[most,breakable]{tcolorbox}

\newtcolorbox{isaframe}[1][]
   { blanker, 
    left=5pt, right=1pt, top=1pt, bottom=1pt,
     borderline west={1pt}{0pt}{cyan},
     before upper=\setlength{\parindent}{0pt},
     fontupper=\ttfamily,
     parbox=true, #1}

\def\A{\mathcal{A}}
\def\wA{\widetilde{\A}}

\def\al{ {\rm alph} }
\def\B{\mathcal{B}}

\def\uu{\mathbf{u}}

\def\N{\mathbb{N}}

\def\uu{\mathbf{u}}

\def\fa{{\rm fact}}

\renewcommand{\restriction}{\mathord{\upharpoonright}}

\newtheorem{theorem}[]{Theorem}
\newtheorem{corollary}[theorem]{Corollary}
\newtheorem{lemma}[theorem]{Lemma}

\newtheorem{definition}[theorem]{Definition}

\theoremstyle{remark}

\crefname{theorem}{Theorem}{The\documentclass{amsart}orems}
\crefname{corollary}{Corollary}{Corollaries}
\crefname{example}{Example}{Examples}
\crefname{lemma}{Lemma}{Lemmas}
\crefname{proposition}{Proposition}{Propositions}
\crefname{definition}{Definition}{Definitions}
\crefname{example}{example}{examples}

\begin{document}
\begin{frontmatter}


\title{The number of primitive words of unbounded exponent in the language of an HD0L-system is finite}



\author[fit]{Karel Klouda}
\author[fit]{Štěpán Starosta}


\affiliation[fit]{organization={Faculty of Information Technology, Czech Technical University in Prague},
            addressline={Thákurova~9}, 
            city={Prague},
            postcode={160 00}, 
            country={Czech Republic}}

\begin{abstract}
Let $H$ be an HD0L-system.
We show that there are only finitely many primitive words $v$ with the property that $v^k$, for all integers $k$, is an element of the factorial language of $H$.
In particular, this result applies to the set of all factors of a morphic word.
We provide a formalized proof in the proof assistant Isabelle/HOL as part of the Combinatorics on Words Formalized project.
\end{abstract}



\end{frontmatter}

\iftoggle{wip}{
\listoftodos
}

\section{Introduction}

Repetitions in languages, or exponents of elements of languages, have been studied in various contexts.
A famous example is the Thue--Morse sequence, which was investigated by A. Thue \cite{Thue12}.
Thue showed that this sequence is overlap-free, meaning it contains no factor of the form $xyxyx$.
By using rational powers of a word, it is possible to write $xyxyx = (xy)^{\frac{|xyxyx|}{|xy|}}$, where $|\cdot|$ denotes the length of a word.
Thus, Thue demonstrated that the Thue--Morse sequence contains no factor of exponent larger than $2$.
This threshold value of $2$ is called the critical exponent.
There have been numerous studies on the topic of critical exponents; see, for instance, \cite{SHALLIT201996} and references therein.

D.~Krieger~\cite{Kri07} presented an algorithm for determining the critical exponent of a fixed point of a non-erasing morphism.
The algorithm operates under the assumption that the critical exponent exists, which is to say that there is a bound on the exponents of factors of the fixed point.
The existence of such a bound is also interesting from the perspective of so-called circularity or recognizability, as discussed in \cite{revekka}. 
Loosely speaking, these two notions describe when one can invert the morphism (on the set of factors of the fixed point).
Specifically, a non-erasing morphism is circular if the exponents in its fixed point are bounded \cite{MiSe, KLOUDA2019131}.
This result applies more generally to D0L-systems, which are instances of L-systems used in \cite{LINDENMAYER1968280} to model growth of organisms.
For further information on this topic also see~\cite{book-of-L}.

In this article, we focus on a generalization of D0L-systems known as HD0L-systems, which were introduced in \cite{Nielsen1974}.
Specifically, we study their factorial language, which is the set of all factors of $\Psi (\varphi^k(w))$ for all integers $k$
with $\varphi$ an endomorphism, $\Psi$ a morpism, and $w$ a word.
These systems add an additional layer to the generation of the language by applying the outer morphism $\Psi$.
In terms of fixed points of morphisms, which are also called purely morphic words, HD0L-systems can be viewed as a generalization of languages of morphic words, i.e., set of factors of an image by $\Psi$ of a fixed point of $\varphi$.

Our main result (\Cref{th:goal}) states that there are finitely many primitive words with unbounded exponent in every such language under the condition that $\Psi$ is non-erasing. 
This result generalizes our previous result of~\cite{KlSt13} which gives an algorithm enumerating the primitive words appearing with unbounded exponent in every language generated by a D0L-system (i.e., for $\Psi$ equal the identity).
Besides this generalization, the present proof works also for $\varphi$ erasing, further generalizing the result of~\cite{KlSt13}.
We also provide a complete formalization, including machine verified proofs, in the proof assistant Isabelle/HOL \cite{Isabelle,IsabelleHOLBook} as part of the Combinatorics on Words Formalized project.

J. Bell and J. Shallit proved \cite{BeSh21} that a factorial language with sublinear factor complexity has finitely many primitive words occurring with unbounded exponent.
In particular, for $k$-automatic words, they gave an algorithm for explicit construction of such words.
This result has a non-empty overlap with our result: the factorial language of an HD0L-system can indeed have sublinear factor complexity, but it can also be larger; for D0L-systems, it follows from J.-J. Pansiot \cite{Pa84} that the factor complexity function can be quadratic, while for HD0L-systems, it follows from R. Deviatov \cite{deviatov}.

The next section contains notation and definitions.
\Cref{sec:L_G} contains results on the language of the underlying D0L-system, which are then used in \Cref{sec:L_H} to prove results on the language of the HD0L-system along with the main theorem.
\Cref{sec:formalization} contains details on the above mentioned formalization of this result, along with the description of its benefits.
\Cref{sec:final} concludes with an open question and a remark on the construction of an algorithm enumerating all primitive words of unbounded exponent.

\section{Preliminaries}

An \textit{alphabet} is a finite set of \textit{letters}.
A~\textit{(finite) word} is a finite sequence over $\A$.
The length of a finite word $w$ is denoted by $|w|$.
The \textit{empty word}, which is the unique word of length $0$, is denoted by $\varepsilon$.
The set of all finite words over the alphabet $\A$ is denoted by $\A^*$.
For $w \in \A^*$, we set $w^0 = \varepsilon$ and $w^{k+1} = w \cdot w^k$ for $k \in \N$, where $\cdot$ is the binary operation of concatenation.
We will often omit the $\cdot$ sign when concatenating two words $w$ a $v$ and simply write $wv$ instead.
A~word $w \in \A^*$ is \textit{primitive} if $v^k = w$ implies $k = 1$.

An \textit{infinite word} is an infinite sequence over $\A$.
A finite word $w$ is a \textit{factor} of a finite or infinite word $z$ if we have $z = pws$ with $p$ being a finite word and $s$ a finite or infinite word.
The word $p$ is a \emph{prefix} of $z$.
If $p$ is non-empty, it is a \emph{strict prefix} of $z$.
Similarly, if $z$ is finite, the word $s$ is its \emph{suffix}.
If $s$ is non-empty, it is a \emph{strict suffix} of $z$.

An infinite word $z$ is \emph{purely periodic} if there exists a word $u$ such that $z = u u u \dots$.
We use the notation $z = u^\omega$.

Given a \emph{language} $L$, i.e., a set of finite words, we say it is \emph{factorial} if it is closed under taking factors:
\[
\forall w \in L, v \text{ is a factor of } w \Rightarrow v \in L.
\]
The smallest superset of $L$ which is factorial is \emph{the factorial closure} of $L$ and is denoted $\fa (L)$.

A \textit{morphism} $\varphi: \A^* \to \B^*$ is a mapping such that for all $v,w \in \A^*$ we have $\varphi(vw) = \varphi(v)\varphi(w)$.
A morphism is \emph{non-erasing} if the only word mapped to the empty word is the empty word.
We set $\| \varphi \| = \max_{a \in \A} |\varphi(a)|$.
If $\A = \B$, we say $\varphi$ is an \textit{endomorphism}.

The triple $G = (\A, \varphi, w_G)$, with $\A$ an alphabet, $\varphi: \A^* \to \A^*$ an endomorphism, and $w_G \in \A^*$, is a \textit{D0L-system}.
The word $w_G$ is usually called the \emph{axiom} of the system.
The language of such system is usually the set $\{\varphi^n(w_G) \mid n\in \N\}$.
However, we are interested in the factorial closure of the set, and hence we define the \emph{factorial language} of $G$ as $L_G = \fa \{\varphi^n(w_G) \mid n\in \N\}$.

A word $w \in \A^*$ is \textit{bounded} if the sequence $|\varphi^n(w)|$ is bounded.
By \textit{bounded letters} we mean elements of $\A$ that form bounded words of length $1$.
If a word or letter is not bounded, we say it is \textit{unbounded}.
If $L_G$ contains infinitely many bounded words, then we say that $G$ is \textit{pushy}; otherwise it is \textit{non-pushy} (see \cite{EhRo83}).

Let $\B$ be an alphabet and $\Psi: \A^* \to \B^*$ a morphism. 
If $G = (\A,\varphi,w_G)$ is a D0L-system, then the 5-tuple $H = (\A,\B,\varphi,\Psi,w_H)$ is an \textit{HD0L-system}.
The word $w_H \in \A^*$ is its \emph{axiom}.
In analogy with the definition of a factorial language of~$G$, we define the \emph{factorial language} of $H$ as $L_H = \fa \left\{ \Psi (\varphi^n(w_H)) \mid n\in \N \right\}$.
Trivially, we have
\[
L_H = \fa (\Psi(L_G)) = \bigcup_{w \in L_G} \fa(\Psi(w)).
\]

The factorial language of a D0L-system is a generalization of a purely morphic language.
A purely morphic language is the set of all factors of an infinite word $\uu$ which is fixed by a morphism $\varphi$, i.e., $\varphi(\uu) = \uu$.
If $a$ is the first letter in $\uu$ and $|\varphi(a)| > 1$, which implies that $|\varphi^n(a)|$ is an unbounded sequence, the set of all factors of $\uu$ coincides with $L_G$ for $G = (\A,\varphi,a)$.
Similarly, a language of an HD0L-system is a generalization of a morphic language; that is, a set of factors of $\Psi(\uu)$, where $\uu$ is fixed by $\varphi$.
Again, if $a$ is the first letter in $\uu$ and $|\varphi(a)| > 1$, the set of all factors of $\Psi(\uu)$ equals $L_H$ for $H = (\A,\B,\varphi,\Psi,a)$.

Given a language $L$, we say that a non-empty word \emph{$w$ is of unbounded exponent in $L$} if $w^k \in L$ for all $k \in \N$.
Since if $w$ is of unbounded exponent in $L$, then so is every power of $w$, it suffices to study primitive words with this property.

In what follows the letter $G$ denotes a D0L-system $(\A, \varphi, w_G)$ and $H$ a HD0L-system $(\A,\B,\varphi,\Psi,w_G)$.
The goal of this article is to study the number of primitive words of unbounded exponent in the factorial language $L_H$.

There is no naive relation between the number of primitive words of unbounded exponent in $L_G$ and $L_H$.
Consider $\varphi: a \mapsto abca, b \mapsto bb, c \mapsto cc$.
Clearly, $b$ and $c$ are primitive words of unbounded exponent in $L_G$ for $G = (\{a,b,c\},\varphi,a)$.
It follows from from~\cite{KlSt13} that these are indeed the only primitive words of unbounded exponent.
Taking $\Psi: a \mapsto 0, b \mapsto 1, c \mapsto 1$, the word $1$ is a primitive word of unbounded exponent in $L_H$ for $H = \left( \{a,b,c\},\{0,1\},\varphi,\Psi,a \right)$.
It can be seen that this is the only such word in $L_H$ and hence, in this case, the number of primitive words of unbounded exponent in $L_H$ is less than the number of the number of primitive words of unbounded exponent in $L_G$.

On the other hand, consider $G$ with the morphism $\varphi: a \mapsto abc, b \mapsto bbc, c \mapsto bac$.
It follows from~\cite{KlSt13} that $L_G$ does not contain any primitive word of unbounded exponent.
Let $\Psi$ be a morphism determined by $a \mapsto 0, b \mapsto 0, c \mapsto 1$.
The word $001$ is a primitive word of unbounded exponent in $L_H$ (in fact, $L_H$ equals the set of all factors of the purely periodic word $001001001\dots$).
That is, in this case, the number of primitive words of unbounded exponent in $L_H$ is greater than the number of primitive words of unbounded exponent in $L_G$.

\subsection{References to the formalization}

In order to allow the reader to refer to our formalization, we give the names of the theorems and auxiliary claims as in the formalization in \texttt{teletype} font.
For instance, the name \texttt{long\_bounded\_factor} is the name of the claim representing \Cref{th:pushy_bounded_factors_format} in our formalization.
Further details on our formalization are described in \Cref{sec:formalization}.

\section{The language $L_G$} \label{sec:L_G}

If $G$ is non-pushy, there exists an upper bound on the length of every bounded word of $L_G$.
If $G$ is pushy, by the definition, no such upper bound exists.
However, from \cite[Proposition 4.7.62]{cassaigne_nicolas_2010}, we know how (sufficiently long) bounded word of $L_G$ looks.
We give here a reformulation of this result following from~\cite[Theorem 12]{KlSt13}.

\begin{theorem}[\texttt{long\_bounded\_factor}] \label{th:pushy_bounded_factors_format}
There exist finite sets $W$ and $U$ such that every bounded $w \in L_G$ is of the form
\begin{equation} \label{eq:long_pushy_form}
w_1u_1^{k_1}w_2u_2^{k_2}w_3   
\end{equation}
where $w_1,w_2,w_3 \in W$, $u_1, u_2 \in U$, and $k_1,k_2 \in \N$.
\end{theorem}

The formulation of \cite[Theorem 12]{KlSt13} is only for pushy $G$, and the words $u_i$ are required to be non-empty.
This requirement forces the statement to be split into 3 cases which are in the reformulation \eqref{eq:long_pushy_form} represented by the following 3 possibilities for pushy $G$: $u_1$ is empty and $u_2$ is non-empty, $u_1$ is non-empty and $u_2$ is empty, or both are non-empty.
The fourth, remaining, case $U = \left\{ \varepsilon \right\}$ covers the case of non-pushy $G$, and since $W$ is finite, is in fact equivalent to it.



We continue with a helpful lemma on $L_G$.
It says that if bounded factors are bounded in length, then a sufficiently long factor always contains a $\varphi^r$-image of an unbounded letter $b$ for a fixed $r$.

\begin{lemma}[\texttt{non\_pushy\_subset\_unbounded\_image}] \label{le:non_pushy_subset_unbouded_image}
Let $r$ be a non-negative integer and $\ell$ an integer.
Let $W$ be a subset of $L_G$.
Assume that for all $v$ and $w$ such that $v$ is a bounded factor of $w \in W$ we have $|v| \leq \ell$.
There exists $n_0$ such that for all $w$ with $w \in W$ and $n_0 \leq |w|$, there exists an unbounded letter $b \in L_G$ such that $\varphi^r (b)$ is a factor of $w$.
\end{lemma}
\begin{proof}
If $W$ is finite, the claim is trivially satisfied.
Assume $W$ is infinite.
Hence the axiom $w_G$ of $G$ is unbounded and $1 < \| \varphi \|$.

Set $n_0 = 2\| \varphi \|^r |w_G| + \ell + 1$.
Assume $w \in W$ and $n_0 \leq |w|$.
Let $t$ be such that $w$ is a factor of $\varphi^t (w_G)$.
The choice of $n_0$ implies that $|w| > \ell$, and thus $w$ is not bounded.
It also implies $r < t$.
Therefore, we have $\varphi^t (w_G) = \varphi^r \varphi^{t-r} (w_G)$.

Let $p'$, $w'$, and $s'$ be such that 
\begin{itemize}
\item $\varphi^r (w') = p' w s'$;
\item $w'$ is a factor of $\varphi^{t - r} (w_G)$;
\item $w'$ is the shortest such factor.
\end{itemize}
The last condition is equivalent to $p'$ being a strict prefix of the $\varphi^r$-image of first letter of $w'$ and $s'$ strict suffix of the $\varphi^r$-image of the last letter of $w'$.
Since $2\| \varphi^r \| \leq 2\| \varphi \|^r < n_0 \leq |w|$ and $\varphi^r (w') = p' w s'$, we conclude that $|w'| > 2$.
Let $w' = yw''x$ with $y$ and $x$ of length $1$ and $w''$ non-empty.
The word $\varphi^r(w'')$ is a factor of $w$.

Now assume that $w''$ is bounded.
It follows that $\varphi^r(w'')$ is also bounded, and we have $|\varphi^r(w'')| \leq \ell$.
Using $\varphi(yw''x) = \varphi^r (w') = p' w s'$, we conclude that $|w| \leq 2\| \varphi^r \| + \ell \leq 2\| \varphi \|^r + \ell < n_0$, which is a contradiction.
Hence, $w''$ is unbounded and there exists an unbounded letter $b$ which is a factor of $w''$.
As $\varphi^r(w'')$ is a factor of $w$, it follows that $\varphi^r(b)$ is a factor of $w$.
\end{proof}

A language $L$ has \emph{uniformly bounded occurrences} of each element if
for all $v \in L$ there exists $k$ such that for all $w \in L$ we have that $|w| \geq k$ implies $v$ is a factor of $w$.

\begin{lemma}[\texttt{bound\_occ\_prim\_unb\_exp}]
If $L_G$ has uniformly bounded occurrences of each element and a primitive word $u$ is of unbounded exponent in $L_H$,
then $L_H$ equals the set of factors of $u^\omega = uuu \dots$.
\end{lemma}
\begin{proof}
Clearly, the set of factors of $u^\omega$ is a subset of $L_H$.

Now assume we have an element $w \in L_H$.
There exists $v \in L_G$ such that $w$ is a factor of $\Psi(v)$.
As the occurrences of $v$ are uniformly bounded in $L_G$, there exists $k$ such that every $x \in L_G$ with $|x| \geq k$ contains $v$ as its factor.

Let $u' = u^{(k+2)\|\Psi\|}$.
Set $z \in L_G$ such that $u'$ is a factor of $\Psi(z)$.
As $u$ is primitive, it is non-empty, and we have $|u'| \geq (k+2)\|\Psi\|$.
Therefore, there exists $z' \in L_G$ such that $|z'| \geq k$ and $\Psi(z')$ is a factor of $u'$.
It follows that $v$ is a factor of $z'$, and thus $\Psi(v)$ is a factor of $u'$, and finally $w$ is a factor of $u'$.
We conclude that $w$ is an element of the set of factors of $u^\omega$.
\end{proof}

In order to understand the language of $L_G$, we define subsets of the alphabet $\A$ that capture the following idea: some letters of $\A$ may appear as factors only at the beginning of the sequence $\left( \varphi^i(w_G) \right)$; those letters are not interesting as they appear only finitely many times; the rest of the letters keep reappearing, and we may furthermore keep track of those subsets that are mapped by a suitable power of $\varphi$ onto words over the same set.
Moreover, only subsets containing unbounded letters are interesting; hence the next definition.

Given a word $w$, we define $\al(w)$ to be the set of letters occurring in $w$.
For a set of words $L$, we define $\al(L) = \bigcup_{w \in L} \al(w)$.

\begin{definition} \label{de:inv_sub}
Let $p$ be a positive integer.
We say that $\A'$ is an \emph{invariant subalphabet with respect to $p$} if $\A' \subset \A$, $\al(\varphi^p(\A')) = \A'$ and $\A'$ contains an unbounded letter.
\end{definition}

In order to find an invariant subalphabet, we first need to fix a positive integer $p$. 
This is done by the next lemma.
For $a \in \A$ and $i \geq 1$, we set $\A_{i,a} = \al(\varphi^i(a))$.


\begin{lemma}[\texttt{pow\_subalph\_l\_ev\_per\_all}] \label{le:p_and_q}
There exist an integer $q$ and a positive integer $p$ such that for every $k, i$ with $k \geq 0$ and $i \geq q$ and every letter $a$, we have 
\[
\A_{i,a} = \A_{i+kp,a}.
\]
Moreover, $p$ can be chosen such that $p \geq q$.
\end{lemma}
\begin{proof}
As $\A$ is finite, the set $\left\{ \A_{i,a}  \mid a \in \A, i \geq 1 \right \}$ is finite.
Hence, given a letter $a$, there exists a positive $p_a$ and an integer $q_a$ such that $\A_{q_a,a} = \A_{q_a+p_a,a}$.
It follows that for every $k \geq 0$ and $i \geq q_a$ we have $\A_{i,a} = \A_{i+kp_a,a}$.
The claim follows for $q = \max_{a \in \A} \left\{ q_a \right\}$ and $p = \prod_{a \in \A} p_a $.
As every multiple of $p$ satisfies the claim, we can choose a suitably large multiple to have $p \geq q$.
\end{proof}

\textbf{Global assumption}: from now on, we assume that $p$ and $q$ are given by~\Cref{le:p_and_q} (with $p \geq q$). 
As $p$ is fixed, we simply use the term \emph{invariant subalphabet} while referring to~\Cref{de:inv_sub}. 
Let $\wA$ be the set of all invariant subalphabets.
If $a$ is unbounded, the alphabet $\A_{p,a}$ is an invariant subalphabet (\texttt{pow\_subalph\_inv\_sub}).
Hence $\wA$ is non-empty, and we can find its minimal elements with respect to the subset relation.
We call these \emph{minimal invariant subalphabets}.

\begin{lemma}[\texttt{inv\_sub\_bot\_pow\_subalph}] \label{le:inv_sub_bot_pow_subalph}
If $\A' \in \wA$ is a minimal invariant subalphabet and $g \in \A'$ is an unbounded letter, then
\[
\A' = \A_{p,g}.
\]
\end{lemma}

\begin{proof}

Since $p \geq q$ and $g$ is an unbounded letter, by~\Cref{le:p_and_q} we have that $\A_{p,g}$ is an invariant subalphabet.
It follows that $\A_{p,g} \subset \A'$.
By the minimality of $\A'$, we conclude that $\A_{p,g} = \A'$.
\end{proof}

For a minimal invariant subalphabet $\A'$ and an unbounded letter $g \in \A'$, consider the triple $G' = (\A', \varphi^p\restriction_{\A'}, g)$.
\Cref{le:inv_sub_bot_pow_subalph} implies $G'$ is indeed a D0L-system.
It forms an interesting subsystem of $G$, which will be helpful to understand $L_G$ and $L_H$.
We first show that all elements of $L_{G'}$ appear with uniformly bounded occurrences for a non-pushy $G'$.

\begin{lemma}[\texttt{inv\_sub\_bot\_bounded\_occ}] \label{le:inv_sub_bot_bounded_occ}
Let $\A' \in \wA$ be a minimal invariant subalphabet and $g \in \A'$ be an unbounded letter.
Set $G' = (\A', \varphi^p\restriction_{\A'}, g)$.
If $G'$ is not pushy, then all elements of $L_{G'}$ appear with uniformly bounded occurrences.
\end{lemma}

\begin{proof}
As $G'$ is not pushy, there exists $\ell_0$ such that all bounded elements of $L_{G'}$ are no longer than $\ell_0$.

By~\Cref{le:non_pushy_subset_unbouded_image} there exists $n_0$ such that for all $w \in L_{G'}$ with $|w| \geq n_0$ there exists an unbounded letter $b$ such that $\varphi^p (b)$ is a factor of $w$.
Since $\A'$ is minimal, we have $b \in \A'$, and hence $g$ occurs in $\varphi^p (b)$.

We show that for all $v \in L_{G'}$ there exists an integer $k$ such that for all $w \in L_{G'}$ with $|w| \geq k$ the word $v$ is a factor of $w$.
Let $v \in L_{G'}$ and let $\ell$ be such that $v$ is a factor of $(\varphi^p)^\ell (g)$.
Set $k = \|\varphi^p\|^\ell \|\varphi\| (n_0+2)$.
Now let $w \in L_{G'}$ be such that $|w| \geq k$.
Let $t$ be such that $w$ is a factor of $(\varphi^p)^t (g)$.
The choice of $k$ implies that $\ell < t$.
Hence, $w$ is a factor of $(\varphi^p)^\ell (\varphi^p)^{t-\ell} (g)$.
Let $p'$, $w'$ and $s'$ be such that
\begin{itemize}
\item $(\varphi^p)^\ell (w') = p' w s'$;
\item $w'$ is a factor of $(\varphi^p)^{t-\ell} (g)$;
\item $w'$ is the shortest such factor.
\end{itemize}
The choice of $k$ implies that $w'$ is of length at least $n_0+2$.
Let $w' = xw''y$ with $x$ and $y$ being of length $1$.
As $|w''| \geq n_0$, it contains the letter $y$.
Since $(\varphi^p)^\ell(w'')$ is a factor of $w$, we conclude that $(\varphi^p)^\ell (g)$ is a factor of $w$, and hence $v$ is a factor of $w$.
\end{proof}

In the next section, we work with sequences of elements of $L_G$.
We distinguish those sequences whose elements contain bounded factors of arbitrary length.

\begin{definition}
We say that a sequence of words $(w_i)$ is \emph{pushy} if for all $n$ there exists $j$ such that $w_j$ contains a bounded factor of length greater than $n$.
\end{definition}

The next lemma states that if we have a non-pushy sequence of elements of $L_G$, with an infinite number of distinct elements, then we can find a minimal invariant subalphabet such that a factor over this subalphabet is longer than some $n$.

\begin{lemma}[\texttt{non\_pushy\_inv\_sub\_bot\_ex}] \label{le:non_pushy_inv_sub_bot_ex}
Let $(w_i)$ be a non-pushy sequence of elements of $L_G$ with $|w_i| \to  +\infty$.
For every $n$, there exists a minimal invariant subalphabet $\A'$, an unbounded letter $g \in \A'$, and integers $j$ and $k$ such that
\[
(\varphi^p)^k (g) \text{ is a factor of } w_j \quad \text{ and } \quad |(\varphi^p)^k (g)| > n.
\]
\end{lemma}

\begin{proof}
Let $k$ be a non-negative integer such that for every unbounded letter $a$ we have
\begin{equation} \label{eq:longer_than_n}
| (\varphi^p)^k (a) | > n.  
\end{equation}
(This claim is formalized as \verb|endomorphism.unbounded_im_all_long|.)

Set $G'(h) = (\A, \varphi^{pk}, \varphi^h (w_G))$.

By the pigeonhole principle, there exists an infinite set of integers $I$ and an integer $h$ such that 
\[
\forall i \in I, w_i \in L_{G'(h)}.
\]\
(This claim is called \verb|endomorphism.pmor_lan_infin_subset_somewhere|.)

As $(w_i)$ is a non-pushy sequence, there exists $\ell$ such that every bounded factor of every $w_i$ is of length at most $\ell$.
We may thus use~\Cref{le:non_pushy_subset_unbouded_image} for $W_I = \left\{ w_i \mid i \in I \right\}$ to obtain $n_0$ such that if $i \in I$ and $|w_i| \geq n_0$, then there exists an unbounded letter $b$ with $b \in L_{G'(h)}$ and
\begin{equation*} 
(\varphi^p)^{k+1} (b) \text{ is a factor of } w_i. 
\end{equation*}

As $I$ is infinite, there exists $j \in I$ with $|w_j| \geq n_0$.
Hence, there exists an unbounded letter $b$ with $b \in L_{G'(h)}$ and
\begin{equation} \label{eq:pe_wj}
(\varphi^p)^{k+1} (b) \text{ is a factor of } w_j.
\end{equation}
The set $\A_{p,b}$ is an invariant subalphabet and it contains a minimal invariant subalphabet $\A'$ as a subset.

Let $g$ be an unbounded letter such that $g \in \A'$.
Then $g$ is a factor of $\varphi^p (b)$.
Hence, from \eqref{eq:pe_wj}, it follows that $(\varphi^p)^k (g)$ is a factor of $w_j$.
Also, by \eqref{eq:longer_than_n}, we have $|(\varphi^p)^k (g)| > n$.
\end{proof}

The last lemma has the following consequence.
As the number of minimal invariant alphabets is finite, and they are themselves finite, 
by the pigeonhole principle, there exists a single minimal invariant subalphabet $\A'$, an unbounded letter $g \in \A'$ and an infinite subset $I$ of integers
such that for every $n \in I$ 
\begin{equation} \label{eq:njk}
\exists j,k  \text{ such that } (\varphi^p)^k (g) \text{ is a factor of } w_j \quad \text{ and } \quad  | (\varphi^p)^k (g) | > n.
\end{equation}
If $n' \not \in I$, we can simply pick some $n \in I$ with $n > n'$ such that \eqref{eq:njk} holds for $n'$.
This is summarized in the following corollary.

\begin{corollary}[\texttt{non\_pushy\_inv\_sub\_bot\_ex\_one}] \label{le:min_inv_sub_rep_ex}
If $(w_i)$ is a non-pushy sequence of elements of $L_G$ with $|w_i| \to  +\infty$,
then there exists a minimal invariant subalphabet $\A'$ and an unbounded letter $g \in \A'$ such that
for all $n$ there exist $j$ and $k$ satisfying the following two conditions:
\[
(\varphi^p)^k (g) \text{ is a factor of } w_j \quad \text{ and } \quad  | (\varphi^p)^k (g) | > n.
\]
\end{corollary}

\section{The language $L_H$} \label{sec:L_H}

Let $u \in L_H$ be of unbounded exponent in $L_H$, i.e., we have $u^k \in L_H$ for every $k$.
The following definition captures the fact that for each $k$, there exists a suitable $\Psi$-preimage of $u^k$ in $L_G$:

\begin{definition}
We say that $(w_i)$ is a \emph{repetition embed sequence} of $u$ if
\begin{itemize}
\item $w_j \in L_G$;
\item $\exists k, \Psi(w_j) \text{ is a factor of } u^k$;
\item $|w_{j+1}| > |w_j|$.
\end{itemize}
\end{definition}

The idea of counting primitive words of unbounded exponent in $L_H$ is based on whether such a word has a pushy repetition embed sequence or not.
The first case is solved by the following lemma.
Recall that a \emph{primitive root} of a non-empty word $w$ is the shortest $v$ such that $w = v^k$ for some $k$.
Clearly, the primitive root of $w$ is primitive.
Two words $u$ and $z$ are \emph{conjugate} if there exists a word $s$ such that $us = sz$.
If $u$ and $z$ are conjugate, we write $u \sim z$.

\begin{lemma}[\texttt{pushy\_rem\_em\_seq\_finite}] \label{le:pushy_embed_finite}
Let $\Psi$ be non-erasing.
The number of primitive words of unbounded exponent in $L_H$ that have a pushy repetition embed sequence is finite.
\end{lemma}
\begin{proof}
By~\Cref{th:pushy_bounded_factors_format}, there exist finite sets $W$ and $U$ such that every bounded $w \in L_G$ is of the form
\begin{equation} \label{eq:long_bounded_factor_applied}
w_1u_1^{k_1}w_2u_2^{k_2}w_3,
\end{equation}
where $w_1,w_2,w_3 \in W$, $u_1, u_2 \in U$ and $k_1,k_2 \in \N$.

Let $\ell_W$ be the maximum length of elements of $W$ and $\ell_U$ the maximum length of elements of $U$.

Assume that a primitive word $u$ is of unbounded exponent in $L_H$ and that $(w_i)$ is a pushy repetition embed sequence of $u$.

Set $k = 2|u|$ and $n = 3 \ell_W + 2k\ell_U$.
As $(w_i)$ is pushy, we may find an index $j$ such that there exists a bounded $w$ which is a factor of $w_j$ and $|w| > n$.
Since $(w_i)$ is a repetition embed sequence, we have that $\Psi(w)$ is a factor of $u^\ell$ for some $\ell$.

As $w$ is bounded, it is of the form \eqref{eq:long_bounded_factor_applied} for some $w_1,w_2,w_3 \in W$, $u_1, u_2 \in U$, and $k_1,k_2 \in \N$.
From this fact combined with $|w| > n = 3 \ell_W + 2k\ell_U$, we conclude that there is $t \in \{1,2\}$ such that $u_t^k$ is non-empty and factor of $w$.

Clearly $\Psi (u_t^k)$ is a factor of $u^\ell$.
As $\Psi$ is non-erasing and $u$ non-empty, we have $|\Psi (u_t^k)| \geq k = 2|u|$.
Hence $\ell \geq 2$.
As $\Psi (u_t^k) = \Psi (u_t)^k$ is a factor of $u^\ell$, it follows that the primitive roots of $\Psi (u_t)$ and $u$ are conjugate, and that the primitive root of $u$ belongs to the set
\[
\left\{ r \mid \exists r', r' \text{ is a primitive root of some element of } \Psi(U) \text{ and } r \sim r' \right\}.
\]
Since the finiteness of $U$ implies finiteness of this set as well, the proof is finished.
\end{proof}

The remaining case is when the primitive word of unbounded exponent in $L_H$ does not have a pushy repetition embed sequence.
By \Cref{le:min_inv_sub_rep_ex}, it implies that the arbitrary power is in fact produced by a subsystem of $G$ determined by minimal invariant subalphabet.
The next lemma states that in such a system, we cannot have an arbitrary number of such primitive words.

\begin{lemma}[\texttt{inv\_sub\_bot\_fin\_prim\_unb\_exp}] \label{le:min_inv_sub_finite}
Let $\Psi$ be non-erasing.
Let $\A' \in \wA$ be a minimal invariant subalphabet.
Let $g \in \A'$ be an unbounded letter and set $G' = (\A', \varphi^p\restriction_{\A'}, g)$.
The number of primitive words of unbounded exponent in $L_{H'} = \fa(\Psi(L_{G'}))$ is finite.
\end{lemma}

\begin{proof}
Assume that $u$ is a primitive word of unbounded exponent in $L_{H'} $ and $(w_i)$ is a repetition embed sequence of $u$.

First assume that $(w_i)$ is not pushy.
Hence, there exists $\ell$ such that all bounded factors of every $w_i$ are not longer than $\ell$.

We show that $G'$ is not pushy.
To obtain a contradiction, assume otherwise, i.e., $G'$ is pushy.
Thus there exists a bounded $w \in L_{G'}$ with $\ell < |w|$.
Let $k$ be an integer such that $w$ is a factor of $(\varphi^p)^k (g)$.
By \Cref{le:non_pushy_subset_unbouded_image} there exists $n_0$ such that for all $j$ with $n_0 \leq |w_j|$ there exists an unbounded letter $b \in L_{G'}$ such that $(\varphi^p)^{k+1} (b)$ is a factor of $w_j$.
As $\A'$ is an invariant subalphabet and $b \in L_{G'}$, we have that $b \in \A'$.
Hence, $g$ is a factor of $\varphi^p(b)$, and consequently, $(\varphi^p)^{k} (g)$ is a factor of $w_j$.
Finally, we conclude that $w$ is a factor of $w_j$, which is a contradiction since no bounded factor $w_j$ is longer than $\ell$. Therefore, $G'$ is not pushy.

Next, we show that every $w \in L_{H'}$ is in fact a factor of $u^s$ for some $s$.
First, let $f$ be an element of $L_{G'}$ such that $w$ is a factor of $\Psi(f)$.
Since $G'$ is not pushy, we may use \Cref{le:inv_sub_bot_bounded_occ} and obtain an integer $t$ such that
\begin{equation} \label{eq:zr}
\forall z \in L_{G'}, t \leq |z|, f \text{ is a factor of } z.
\end{equation}
As $(w_i)$ is a repetition embed sequence of $u$, we may find $i$ such that $t \leq |w_i|$ and $w_i$ is a factor of $\varphi^{p(d+1)} (g)$ for some $d$.
By \eqref{eq:zr}, the word $f$ is a factor of $w_i$.
Since $\Psi (w_i)$ is a factor of $u^s$ for some $s$, we conclude that $\Psi(f)$ is a factor of $u^s$.
As $w$ is a factor of $\Psi(f)$, the word $w$ is also a factor of $u^s$.

Let $u'$ be a primitive word of unbounded exponent in $L_{H'}$.
Let $k'$ be such that $|u'^{k'}| \geq |u^2|$ and $k' \geq 2$.
As $u'^{k'}\in L_{H'}$, it is a factor of $u^s$ for some $s$.
As both $u$ and $u'$ are primitive, it follows that they are conjugate.
We conclude that the number, up to conjugation, of primitive words of unbounded exponent of $L_{H'}$ that have a non-pushy repetition embed sequence is 1.

Since by \Cref{le:pushy_embed_finite}, the number of repetitions of $L_{H'}$ that have a pushy repetition embed sequence is also finite, the proof is finished.
\end{proof}

It remains to put the two cases together.

\begin{theorem}[\texttt{mor\_lan\_fin\_prim\_unb\_exp}] \label{th:goal}
Let $\Psi$ be non-erasing.
The number of primitive words of unbounded exponent in $L_H$ is finite.
\end{theorem}
\begin{proof}
Let $u$ be a primitive word of unbounded exponent in $L_H$.
Let $(w_i)$ be a repetition embed sequence of $u$.

If $(w_i)$ is pushy, then by~\Cref{le:pushy_embed_finite}, the number of such factors $u$ is finite.
If $(w_i)$ is not pushy, then it follows from~\Cref{le:min_inv_sub_rep_ex} that
there exists a minimal invariant subalphabet $\A' \in \wA$ and unbounded letter $g \in \A'$ such that $u$ is a primitive element of unbounded exponent in $G' = (\A', \varphi^p\restriction_{\A'}, g)$.
As the number of primitive elements of unbounded exponent in such system $G'$ is finite by~\Cref{le:min_inv_sub_finite}, and there is a finite number of such subsystems, the number of such factors $u$ is also finite.
\end{proof}

\section{Notes on the formalization} \label{sec:formalization}

In this section, we first briefly summarize the benefits of mathematics formalization in general.
We give description of our formalization and we conclude with specific achievements of our formalization.

The goals of mathematics formalization can be perceived in two steps.
First, the goal is to write down definitions, statements, and proofs in a formal language
with the proofs using a fixed set of inference rules.
The second goal is to verify the correctness of the formalized proofs algorithmically.

The benefits of formalization follow directly from achieving these goals.
The usage of a formal language removes ambiguity, and thus results in consistency.
It also provides better presentation and higher reusability.
The reusability stems also from the simplicity and generality of definitions and claims, and from their organization;
as reusability becomes a requirement in a larger formalization, more general definitions and claims, and better structure may be considered as other benefit of a formalization.

The ability of having a proof checked algorithmically, possibly by a machine, is the most obvious benefit; especially for proofs which are long and technical---beyond what a human reader is capable of checking.
The correctness of a proof is then given by the correctness of the checking system.
For a more detailed overview of benefits of mathematics formalization see \cite{why_formalize}.

Our formalization is part of a larger formalization effort, specifically the Combinatorics on Words Formalized Project \cite{CoW_gitlab}.
It is done in the proof assistant Isabelle/HOL \cite{Isabelle,IsabelleHOLBook}.
An overview of the fundamental part of the project, the library containing tools to work with finite words, is introduced in \cite{itp2021}.

The currently described formalization is archived at~\cite{zenodo_1_10_1}.
A most up-to-date version is available at the project's public repository \cite{CoW_gitlab}.
It is part of the session called ``CoW\_Infinite''.
A session is a collection of basic formalization units, called theories.
The presented result is in the theory called ``Morphic\_Language\_Unbounded\_Exponent''.
The general results on factorial languages of HD0L-systems rely on ``Languages'' theory, described in more detail in~\cite{starosta2023infinite}.
Factorial languages of HD0L-systems are formalized as a locale called ``morphic\_language''.
A \emph{locale} is Isabelle's mechanism to avoid repetition of assumptions, i.e., it allows to fix them and formalize numerous claims under these assumption.
It also allows easy reuse of these facts; see the formalized version of \Cref{th:goal} below.

Our formalization is sometimes slightly more general than what is written above.
As already mentioned above, this is one of the benefits of formalization; one tries for the formalized facts to be reusable later which includes simplifying the assumptions as much as possible and generalizing the claims.
The top example of this effect is that the assumption on $\varphi$ to be non-erasing was dropped in comparison to what we had before the formalization process started.
Of course, the price for dropping this assumption is the need for more complicated proofs.
For instance, the claim producing \eqref{eq:longer_than_n} above is trivial for non-erasing $\varphi$.

There are also minor differences in the notation: in the formalization, the two morphisms are \texttt{f} and \texttt{h} rather than $\varphi$ and $\Psi$.

\Cref{th:goal} is called \texttt{morphic{\isacharunderscore}{\kern0pt}language{\isacharunderscore}{\kern0pt}finite{\isacharunderscore}{\kern0pt}unbounded{\isacharunderscore}{\kern0pt}exponent} in the formalization  and reads as follows:

\begin{isaframe}
\isacommand{theorem}\isamarkupfalse%
\ morphic{\isacharunderscore}{\kern0pt}language{\isacharunderscore}{\kern0pt}finite{\isacharunderscore}{\kern0pt}unbounded{\isacharunderscore}{\kern0pt}exponent{\isacharcolon}{\kern0pt}\isanewline
\ \ \isakeyword{assumes}\ {\isachardoublequoteopen}finite\ {\isacharparenleft}{\kern0pt}UNIV{\isacharcolon}{\kern0pt}{\isacharcolon}{\kern0pt}{\isacharprime}{\kern0pt}a\ set{\isacharparenright}{\kern0pt}{\isachardoublequoteclose}\ \isakeyword{and}\isanewline
\ \ {\isachardoublequoteopen}endomorphism\ f{\isachardoublequoteclose}\ {\isachardoublequoteopen}nonerasing{\isacharunderscore}{\kern0pt}morphism\ h{\isachardoublequoteclose} \isanewline
\ \ \isakeyword{shows}\ {\isachardoublequoteopen}finite\ {\isacharbraceleft}{\kern0pt}u{\isachardot}{\kern0pt}\ primitive\ u\ {\isasymand}\isanewline  {\isacharparenleft}{\kern0pt}{\isasymforall}n{\isachardot}{\kern0pt}\ u\ \isactrlsup {\isacharat}{\kern0pt}\ n\ {\isasymin}\ morphic{\isacharunderscore}{\kern0pt}language\ f\ h\ {\isacharparenleft}{\kern0pt}axiom{\isacharcolon}{\kern0pt}{\isacharcolon}{\kern0pt}{\isacharprime}{\kern0pt}a\ list{\isacharparenright}{\kern0pt}{\isacharparenright}{\kern0pt}{\isacharbraceright}{\kern0pt}{\isachardoublequoteclose}
\end{isaframe}

The assumption \texttt{finite {\isacharparenleft}{\kern0pt}UNIV{\isacharcolon}{\kern0pt}{\isacharcolon}{\kern0pt}{\isacharprime}{\kern0pt}a\ set{\isacharparenright}
} represents the assumption on finiteness of the alphabet; the alphabet is given by the type variable \texttt{{\isacharprime}{\kern0pt}a} for letters, and the assumption says that the universe for this type variable is finite.
The expression \texttt{u\ \isactrlsup {\isacharat}{\kern0pt}\ n} is the notation for the \texttt{n}-the power of the word 
\texttt{u}, i.e., the power $u^n$ as defined as above.
Hence, the claim following the keyword \isakeyword{shows} can be read as ``the set of words \texttt{u} such that \texttt{u} is primitive and for all \texttt{n}, the \texttt{n}-th power of the word \texttt{u} belongs to 
\texttt{morphic{\isacharunderscore}{\kern0pt}language\ f\ h\ {\isacharparenleft}{\kern0pt}axiom{\isacharcolon}{\kern0pt}{\isacharcolon}{\kern0pt}{\isacharprime}{\kern0pt}a\ list{\isacharparenright}}, which is 
the factorial language of the HD0L-system $(\A,\B,\texttt{f},\texttt{h},\texttt{axiom})$ with the alphabets $\A$ and $\B$ being given by the types of \texttt{f} and \texttt{h}.

Our formalization is a witness of the above mentioned benefits of mathematics formalization.
Besides the mentioned more general claims, the process of formalization resulted in many reusable auxiliary claims and in an almost complete reorganization of the proof, ending in a better presentation of the written article.
Last but not least, the proofs provided above are verified and verifiable by a machine.

\section{Final remarks} \label{sec:final}

\subsection{The case of $\Psi$ erasing}
The question on whether the assumption of $\Psi$ being non-erasing can be dropped is very natural.
However, we are not able to provide a proof using our current techniques, which rely on finding a subalphabet such that we are sure to have a power with any exponent over this subalphabet.
We have two types of subalphabets: one is produced by bounded letters, and it is covered by \Cref{le:pushy_embed_finite};
the other is a minimal invariant subalphabet, covered by \Cref{le:min_inv_sub_finite} (which reuses \Cref{le:pushy_embed_finite} in its proof).
Allowing $\Psi$ to be erasing implies that we are no longer able to find this power of any length in a specific subalphabet.
In other words, we are not able to exclude the case where a power with any exponent never occurs inside one of these subalphabets as this fact is no longer projected by $\Psi$ from $L_G$.

Despite the fact that we failed to lift this assumption, we conjecture it can indeed be dropped.

\subsection{An algorithm to count primitive words of unbounded exponents of $L_H$} 

With a few more ingredients, the proof of~\Cref{th:goal} is in fact constructive and we can retrieve an algorithm enumerating all primitive words of unbounded exponent of $L_H$ up to conjugation.
As already recalled above, there are two types of such primitive words.
The two cases are in fact distinguished by the pushiness of their repetition embed sequences.

To enumerate the words that have a pushy embed sequence, we need to determine the set $U$ of \Cref{th:pushy_bounded_factors_format}, and then find its elements that occur in $L_G$ with unbounded exponent.
This is part of the algorithm of \cite{KlSt13} finding primitive words that appear in $L_G$ with unbounded exponent.
As already mentioned, the algorithm requires $\varphi$ to be non-erasing, so its refinement to a general $\varphi$ is needed.
Following the proof of~\Cref{le:pushy_embed_finite}, it suffices to determine the primitive roots of $\Psi(U)$ and keep those that are unique up to conjugation.
Since for every $u \in U$ there is $k$ such that $\varphi^k(u) = u$, it can be shown that the number of elements of $U$ is at most the maximum largest common multiple of each partition of the number of bounded letters in $\mathcal A$.
This bound is given by Landau's function (see e.g.,~\cite{Massias1984}).

To enumerate the words that do not have a pushy embed sequence, we need to enumerate all minimal invariant subalphabets such that the languages of the systems they induce are subsets of $L_H$.
This is equivalent to checking whether for a given minimal $\A' \in \wA$ there is an unbounded letter $b \in \A'$ that occurs in $L_G$.
As it follows from the proof of~\Cref{le:min_inv_sub_finite}, each such minimal invariant subalphabet can produce at most one primitive word of unbounded exponent containing an unbounded letter up to conjugation.

\section*{Acknowledgements}

K. Klouda was supported by The Ministry of Education, Youth and Sports of the Czech Republic through the project\linebreak CZ.02.1.01/0.0/0.0/16\_019/0000778.
Š. Starosta acknowledges support by the Czech Science Foundation grant GA\v CR 20-20621S.

\bibliographystyle{siam}
\bibliography{biblio}

\end{document}